\documentclass[
12pt,
a4paper]{article}

\usepackage{color}

\usepackage[english]{babel}
\usepackage{amsmath}
\usepackage{amssymb}

\usepackage{amsfonts}

\usepackage{array}
\usepackage[dvips]{graphicx}

\usepackage{epsfig}
\usepackage{float}

\usepackage{theorem}

\font\msbmten=msbm10 scaled 1100
\def\Bbt#1{\msbmten#1}

\newcommand{\PP}{\hbox{\Bbt P}}

\newcommand{\calA}{{\mathcal A}}

\newcommand{\calF}{{\mathcal F}}

\DeclareMathOperator{\bE}{\Bbt{E}}

\newtheorem{proposition}{Proposition}[section]

\newtheorem{remark}{Remark}[section]
\newtheorem{lemma}{Lemma}[section]
\newtheorem{theorem}{Theorem}[section]

\newenvironment{pf}[1][Proof]{\noindent{\bf #1\/}
\begin{small}
\noindent%
}{\hspace*{\fill}\nolinebreak[1]%
\hspace*{\fill}\EndProofMarker\\\medskip
\end{small}}
\newcommand{\EndProofMarker}{$\Box$}

\numberwithin{equation}{section}

\begin{document}
\title{A Central Limit Theorem,\\ and related results,\\
for a two-color randomly reinforced urn}
\author{Giacomo Aletti$^\star$, Caterina May$^\natural$ and Piercesare Secchi$^\sharp$}
\date{\today}
\maketitle
 \begin{center}
 {\small
 $^\star$Dipartimento di Matematica ``F. Enrigues''\\
   Universit\`a degli Studi di Milano\\
   via Cesare Saldini 50, 20133 Milano, Italy\\
   {\tt giacomo.aletti@mat.unimi.it}\\
   \medskip

 $^\natural$Dipartimento SEMEQ\\
   Universit\`a del Piemonte Orientale\\
   via Perrone 18, 28100 Novara, Italy\\
   {\tt caterina.may@eco.unipmn.it}\\
   \medskip

 $^\sharp$ MOX--Dipartimento di Matematica \\
         Politecnico di Milano \\
         via Bonardi 9, 20133 Milano, Italy\\
 {\tt piercesare.secchi@polimi.it}}
 \end{center}

 \noindent {\bf Keywords}: reinforced processes, generalized Polya urns, convergence of conditional distributions.
\vspace*{0.5cm}

 \noindent {\bf 2000 AMS Subject Classification}:
60F05
 \vspace*{0.5cm}

\begin{abstract} We prove a Central Limit Theorem for the sequence of random compositions of a two-color randomly reinforced urn. As a consequence, we are able to show that the distribution of the urn limit composition has no point masses.
\end{abstract}

\section{Introduction}
Consider an urn containing initially \(x\) balls of color black and \(y\) balls of color white, with \(x\) and \(y\) non negative real numbers such that \(x+y>0.\) The urn is sequentially sampled: whenever the color of the sampled ball is black, the ball is replaced in the urn together with a random number of black balls, generated at that instant from a distribution \(\mu\) with non negative bounded support; whenever the sampled ball is white, the ball is replaced in the urn together with a random number of balls, generated at that instant from a distribution \(\nu\) with non negative bounded support. This is an informal description
of the Randomly Reinforced Urn (RRU) introduced in \cite{Muliere.et.al.06} and studied  in \cite{AlettiMaySecchi07,Beggs05,DurhamYu,Hopkins.Posch05,LiDuhramFlournoy,MayFl07,May.et.al.05}  under various assumptions concerning the reinforcement distributions \(\mu\) and \(\nu.\) The urn has an interesting potential for applications since it describes a general model for reinforcement learning (\cite{Beggs05,Hopkins.Posch05}); in clinical trials, it implements an optimal response adaptive design (\cite{DurhamFlournoyLi,MayFl07,Sco07,PS07}).

The focus of this paper is on the asymptotic behavior of the sequence \(\{Z_n\}\) describing the random proportions of black balls in the urn along the sampling sequence; in \cite{Muliere.et.al.06} it is proved that the sequence \(\{Z_n\}\) converges almost surely to a random limit \(Z_\infty \in [0,1].\)

When \(\mu=\nu,\) a RRU is a special case of the generalized Polya urn studied by Crimaldi in \cite{Crimaldi08}; for the sequence of random proportions \(\{Z_n\}\) generated by her urn, Crimaldi proves a Central Limit Theorem by showing almost sure conditional convergence to a Gaussian kernel of the sequence \(\{\sqrt{n}(Z_n-Z_\infty)\}.\) Crimaldi's result does not hold for a general RRU; in this paper we extend it to cover the case of a RRU with reinforcement distributions \(\mu\) and \(\nu\) having the same mean. When the means of \(\mu\) and \(\nu\) are different, the limit proportion \(Z_\infty\) of a RRU is a point mass either in 1 or in 0, according to the reinforcement distribution having the larger mean, as proved with different arguments in \cite{Beggs05,Hopkins.Posch05,Muliere.et.al.06}.

A nice implication of our RRU Central Limit Theorem is that we are now able to prove that the distribution of the limit proportion \(Z_\infty\) has no point masses in \([0,1],\) when the means of the reinforcement distributions are the same. This gives a new drive to the problem concerning the absolute continuity of the distribution of the limit proportion of a generalized Polya urn, considered, for instance, also in \cite{Pemantle90}.

The paper is organized as follows. In the next section we will formally introduce the RRU model along with the notation used in the paper. The main results of the paper are stated in section 3, while proofs appear in the following section.
A remark on the absolute continuity of the distribution of \(Z_\infty\) concludes the paper.

\section{Model description and notations}
On a rich enough probability space \((\Omega, \cal A, P),\) define two independent infinite sequences of random elements,
\(\{U_n\}\) and \(\{(V_n,W_n)\}\); \(\{U_n\}\) is a sequence of i.i.d. random variables uniformly distributed on \([0,1],\) while
\(\{(V_n,W_n)\}\) is a sequence of i.i.d bivariate random vectors with components uniformly distributed on \([0,1].\)
Given two probability distributions $\mu$ and $\nu$ on \([0,\beta],\) with \(\beta>0,\) indicate their quantile functions with
$q_\mu$ and $q_\nu,$ respectively. Then, define an infinite sequence \(\{(R_X(n),R_Y(n))\}\) of bivariate random vectors by setting,
for all \(n,\)
\[
R_X(n)=q_\mu(V_n)\;\;\mbox{and}\;\; R_Y(n)=q_\nu(W_n).
\]
Note that, whereas the sequences \(\{(R_X(n),R_Y(n))\}\) and \(\{U_n\}\) are independent, the random variables \(R_X(n)\) and \(R_Y(n)\) might be dependent; however, for every \(n,\) their distributions are \(\mu\) and \(\nu,\) respectively. We indicate with \(m_\mu\) and \(m_\nu,\) and with \(\sigma^2_\mu\) and \(\sigma^2_\nu,\) the means, and the variances, of two random variables \(R_X\) and \(R_Y\) having probability distributions \(\mu\) and \(\nu,\) respectively.

We are now ready to introduce a process whose law is that of a Randomly Reinforced Urn as defined in \cite{Muliere.et.al.06} .
Let \(x\) and \(y\) be two non-negative real numbers such that
\(x+y>0.\)
Set  $X_0=x$, $Y_0=y$, and, for $n= 0,1,2,...$, let
\begin{equation}\label{def1}
\left \{ \begin{array}{lll}
  X_{n+1}&=&X_n+R_{X}(n+1)\delta_{n+1},\\
  Y_{n+1}&=&Y_n+R_{Y}(n+1)(1-\delta_{n+1}),
\end{array} \right.
\end{equation}
where the variable \(\delta_{n+1}\) is the indicator
of the event \(\{U_{n+1}\leq X_{n}(X_{n}+Y_{n})^{-1}\}.\)
The law of \(\{(X_n,Y_n)\}\) is that of the
stochastic process counting, along the sampling sequence, the number of black and
white balls present in a RRU with initial
composition \((x,y)\) and reinforcement distributions equal
to \(\mu\) and $\nu,$ respectively.

For \(n=0,1,2,\ldots\) let
$$Z_n=\dfrac{X_n}{X_n+Y_n};$$ \(Z_n\) represents the proportion
of black balls in the urn before the \((n+1)\)-th ball is sampled from it.
In \cite{Muliere.et.al.06} it is proved that \(\{Z_n\}\) is eventually a bounded sub- or super-martingale, according to the mean of \(\mu\) being larger or smaller than that of \(\nu.\) Hence, for \(n\) growing to infinity, \(Z_n\) converges almost surely, and in $L^p$, $1\leq p\leq\infty,$ to a random variable \(Z_\infty \in [0,1].\)

For \(n=1,2,...,\) let $R_n= \delta_n R_X(n)+(1-\delta_n)R_Y(n)$ be the urn reinforcement, when the urn is sampled for the \(n\)-th time, and set
\[
\begin{gathered}
Q_{n-1}^X = \frac{R_X(n)}{\sum_{i=1}^{n} R_i},
\qquad Q_{n-1}^Y = \frac{R_Y(n)}{\sum_{i=1}^{n} R_i},\\
Q_{n-1} = \frac{R_n}{\sum_{i=1}^{n} R_i} = \delta_n Q_{n-1}^X+(1-\delta_n)Q_{n-1}^Y,
\end{gathered}
\]
with $Q_{n-1}^{X}=Q_{n-1}^{Y}=Q_{n-1}=1$ if $R_i=0$ for all $i=1, \dots, n$.
For shortness, we will write \(D_n\) for the random number $X_n+Y_n,$
interpreted as the size of the urn before it is sampled for the \((n+1)\)-th time.
Clearly, $D_0=x+y$ while $D_{n+1}=D_n+R_{n+1}$, for \(n=0,1,2,....\)
Finally let $\calA_n =\sigma(U_1,\ldots,U_n,(V_1,W_1),\ldots,(V_n,W_n))$ and consider the filtration \(\{\calA_n\};\)
for \(n=1,2,...,\) we indicate with $M_n$ and $A_n$ the two terms given by the Doob's semi-martingale decomposition
of $Z_n$: i.e.
\[
Z_n= Z_0 + M_n+ A_n,
\]
where $\{M_n\}$ is a zero mean martingale with respect to \(\{\calA_n\},\) while $\{A_n\}$ is previsible with respect to \(\{\calA_n\}.\) Theorem 2 in  \cite{Muliere.et.al.06} shows that $\{A_n\}$ is eventually increasing or decreasing.

\section{Main results}
For every set \(A \in \cal A,\) every \(\omega \in \Omega\) and \(n=1,2,...,\)  define
\[
K_n(\omega, A)= P(\sqrt{n}(Z_n-Z_\infty)\in A|\calA_n)(\omega);
\]
i.e. \(K_n\) is a version of the conditional distribution of \(\sqrt{n}(Z_n-Z_\infty)\) given \(\calA_n.\)
When the reinforcement distributions of an RRU are the same, i.e. \(\mu=\nu,\) and \(\mu\) is different from the point mass at 0, Corollary 4.1 in
\cite{Crimaldi08} shows that, for almost every \(\omega \in \Omega,\)
the sequence of probability distribution \(\{K_n(\omega, \cdot)\}\) converges weakly to
the Gaussian distribution \[N(0, hZ_\infty(\omega)(1-Z_\infty(\omega))),\] where
\[
h=\frac{\int_0^\beta k^2 \mu(dk)}{(\int_0^\beta k \mu(dk))^2}.
\]
The next theorem extends this Central Limit result to a general RRU with
reinforcement distributions \(\mu\) and \(\nu\) having the same mean.

\begin{theorem}\label{thm:mainExtension}
Assume $m_\mu = m_\nu =m >0.$
Let
$$
  H = m^{-2}\left(Z_{\infty}\int_0^\beta k^2 \nu(dk)+(1-Z_{\infty})\int_0^\beta k^2 \mu(dk)\right).
$$
Then, for almost every $\omega\in\Omega$,
the sequence  of probability distributions $\{K_n(\omega,\cdot)\}$ converges
weakly to the Gaussian distribution
$$ N \big(0,H(\omega)Z_\infty(\omega)(1-Z_\infty(\omega))\big).$$
\end{theorem}

When \(\mu=\nu,\) the distribution of \(Z_\infty\) has no point masses; this has been proved in \cite{May.et.al.05}.
May and Flournoy show in \cite{MayFl07} that equality of the means of the reinforcement distributions is a sufficient condition for proving that
$\PP(Z_\infty=\{0\})=\PP(Z_\infty=\{1\})=0.$ As a nice consequence of Theorem \ref{thm:mainExtension}, we are now able to prove that the distribution of \(Z_\infty\) has no point masses, when the means of \(\mu\) and \(\nu\) are the same.

\begin{theorem}\label{theo:nomasses}
If $m_\mu = m_\nu=m>0,$ the distribution of $Z_\infty$ has no point masses. I.e.
$\PP(Z_\infty=\{x\})=0,$ for all $x\in [0,1].$
\end{theorem}

When $m_\mu \not=m_\nu,$
the distribution of \(Z_\infty\) is the point mass at 1 or at 0,
according to whether \(m_\mu\) is larger or smaller than \(m_\nu\); this has been proved
in \cite{Beggs05,Hopkins.Posch05,Muliere.et.al.06} under the assumption that the supports of \(\mu\) and \(\nu\) are bounded away from 0.
Within the framework of the present paper, we are able to show that the result holds more generally when the supports of \(\mu\) and \(\nu\) are contained in the interval \([0,\beta].\)

\begin{theorem}\label{cor:extension}
Assume $ m_\mu > m_\nu.$ Then $\PP(Z_\infty=1)=1.$
\end{theorem}

\section{Proofs and auxiliary results}
The proof of Theorem \ref{thm:mainExtension} will make use of a few auxiliary results, that we state and prove as Lemmas
\ref{cor:defiBino1} - \ref{lem:conv_sum}.

From now on, given a sequence $\{\xi_n\}$ of random variables, we will denote by $\Delta \xi_{n+1}$ the increments $(\xi_{n+1}-\xi_{n})$. Moreover, given any two sequences $\{a_n\}$ and $\{b_n\}$ of real numbers, we will use the symbol $a_n \backsim b_n$ to denote that $a_n/b_n \rightarrow 1 $ as $n \rightarrow \infty.$

\begin{lemma}\label{cor:defiBino1}
Let \(R\) and \(D\) be  two random variables defined on \((\Omega, \calA, P)\) with values in \(B_R=[0,\beta]\) and \(B_D=[0,\infty),\) respectively, and let \({\cal G}\) be a sub-sigma-field of \(\calA\) such that \(R\) is independent of \({\cal G}\) while \(D\) is measurable with respect to \({\cal G}.\)
Let \(h\) be a measurable real valued function defined on \(B_R\times B_D\) and such that \(h(\cdot,t)\) is convex for all \(t \in B_D.\)
Then, for almost every \(\omega \in \Omega,\)
\[
h(\bE(R),D(\omega)) \leq \bE \big( h(R,D) \big| {\cal G} \big)(\omega) \leq
\Big(\frac{\bE(R)}{\beta} h(\beta,D(\omega)) +
\frac{\beta-\bE(R)}{\beta} h(0,D(\omega))\Big) .
\]
The previous inequalities are reversed if $h(\cdot,t)$
is concave for all \(t \in B_D.\)
\end{lemma}

\begin{pf}
If \(\pi\) is the probability distribution of \(R,\)
\[
\bE \big( h(R,D) \big| {\cal G} \big)(\omega) = \int_0^\beta h(x, D(\omega)) \pi(dx)
\]
for almost every \(\omega \in \Omega.\)
The left inequality is now an instance of Jensen's inequality. The right inequality follows
after noticing that
\[
h(x,t) \leq x\frac{h(\beta,t)}{\beta} + \frac{\beta-x}{\beta}h(0,t)
\]
for all \((x,t) \in B_R \times B_D,\) since \(h(\cdot,t)\) is convex.
\end{pf}

As a consequence of the previous Lemma, we  can bound the increments
$\Delta A_n$ of the compensator process \(\{A_n\}\).
First note that, for all \(n=0,1,2,...,\)
\begin{equation*}\label{eq:eqDeltaAn_1}
\Delta {A}_{n+1}
=
\bE(\Delta Z_{n+1}|{{\mathcal A}_n}) = Z_{n}(1-Z_n)
A^*_{n+1}
\end{equation*}
almost surely, where
\begin{equation}\label{eq:A^*}
A^*_{n+1}=\bE\Big( \frac{\frac{R_X(n+1)}{D_n}}{1+\frac{R_X(n+1)}{D_n}} -
\frac{\frac{R_Y(n+1)}{D_n}}{1+\frac{R_Y(n+1)}{D_n}} \Big|{{\mathcal A}_n}
\Big).
\end{equation}

\begin{lemma}\label{lem:A^*}
If  $m_\mu > m_\nu$,  for almost every \(\omega \in \Omega\) there is an $a>0$ such that
\begin{equation*}
A^*_{n+1}(\omega) \geq \frac{a}{D_n(\omega)}
\end{equation*}
eventually.
If $m_\mu=m_\nu=m>0$, for \(n=0,1,2...,\) and almost every \(\omega \in \Omega,\)
\begin{equation*}
|A^*_{n+1}(\omega)| \leq
\frac{m(\beta-m)}{(\beta+D_n(\omega))(m+D_n(\omega))}=O \left(\frac{1}{D_n^2(\omega)}\right).
\end{equation*}
\end{lemma}
\begin{pf}
Note that
\(
h(x,t)=
\frac{x}{x+t},
\) is a concave bounded function
of $x\geq 0$ for any fixed $t\geq 0$.
For \(n=0,1,2,...,\)
\[
A^*_{n+1} = \bE\big( h(R_{X}(n+1),D_n) -h(R_{Y}(n+1),D_n) \big|{{\mathcal A}_n}
\big);
\]
hence, by applying Lemma~\ref{cor:defiBino1} we get
\begin{align*}
& A^*_{n+1}\geq \frac{m_\mu}{\beta+D_n} - \frac{m_\nu}{m_\nu+D_n}
=\frac{D_n(m_\mu-m_\nu)-m_\nu(\beta-m_\mu)}{(m_\nu+D_n)(\beta+D_n)}
\\
\intertext{and}
\\
& A^*_{n+1}\leq \frac{m_\mu}{m_\mu+D_n} - \frac{m_\nu}{\beta+D_n}
=\frac{D_n(m_\mu-m_\nu)+m_\mu(\beta-m_\nu)}{(m_\mu+D_n)(\beta+D_n)},
\end{align*}
on a set of probability one.
The thesis is now a consequence of the fact that $\lim_{n \rightarrow \infty} D_n=\infty$ almost surely
(see, e.g., \cite[Proposition~2.3, Proposition~2.4]{MayFl07}).%
\end{pf}

Indeed, when \(m_\mu=m_\nu=m>0,\) two educational cases emerge by inspection of \(A^*_{n+1}\) in the light of Lemma \ref{cor:defiBino1}.
The first one is when \(\mu\) is the point mass at \(m:\) then \(A^*_{n+1} \geq 0\) for all \(n=0,1,2,...\) and the process \(\{Z_n\}\) is a bounded submartingale. On the other extreme, let \(\mu\) be the distribution of the random variable \(\beta \zeta,\) with \(\zeta\) distributed according to a Bernoulli\((m/\beta);\) then \(A^*_{n+1} \leq 0\) for all \(n=0,1,2,...\) and the process \(\{Z_n\}\) is a bounded supermartingale.

In \cite[Lemma~A.1(iii)]{MayFl07} it is proved that \(\lim_{n \rightarrow \infty} D_n/n=m\) almost surely, when \(m_\mu=m_\nu=m>0.\) The next Lemma improves our general understanding of the growth speed of the urn size \(D_n.\)

\begin{lemma}\label{lem:MayGarciaBis}
Assume that \(\min(m_\mu,m_\nu)>0.\)
For all $c,\alpha\geq 0$, there are two constants $0<a_1<a_2<\infty$ such that
\begin{equation}\label{eq: d_asympt_bounds}
\frac{a_1}{n^\alpha} \le \bE\frac{1}{(c+D_n)^\alpha} \le \frac{a_2}{n^\alpha}
\end{equation}
eventually.
Moreover, if $m_\mu=m_\nu=m>0$, then
\[
\bE\frac{1}{(c+D_n)^\alpha} \sim \frac{1}{(c+D_0+mn)^\alpha}.
\]
\end{lemma}

\begin{pf}
It is trivial to prove the lemma when the supports of \(\mu\) and \(\nu\) are both bounded away from 0; this is the case, for instance, when \(\mu\) and \(\nu\) are both point masses at real numbers different from 0.
For the general case
set $\sigma^2=\min(\sigma^2_\mu,\sigma^2_\nu)$ and assume that $\max(\sigma^2_\mu,\sigma^2_\nu)>0;$ without loss of generality we also assume that \(m_\mu\geq m_\nu>0.\)

The left inequality in (\ref{eq: d_asympt_bounds}) follows from Jensen's inequality:
\[
 \bE\frac{1}{(c+D_n)^\alpha}\geq \frac{1}{(c+\bE (D_n))^\alpha}
\geq \frac{1}{(c+D_0 +nm_\mu)^\alpha}.
\]

For proving the right inequality, we consider two cases.
\medskip

\noindent
\emph{Case 1: $\sigma>0$}. For $n=1,2,...$ and $i=1,\ldots,n$ set
\begin{equation}\label{eq:mart_diff_arr}
L_{ni}= \frac{\delta_i(R_X(i)-m_\mu)+
(1-\delta_i)(R_Y(i)-m_\nu)}{\sqrt{n}
\sqrt{Z_{i-1}\sigma^2_\mu+(1-Z_{i-1})\sigma^2_\nu}}
\leq \frac{R_i-m_\nu}{\sqrt{n}\sigma}.
\end{equation}
Then one can show that
$\{L_{ni},\calF_{ni}=\calA_{i},n=1,2,..., i=1,\ldots,n\}$ is a martingale
difference array such that, for all \(n=1,2,...,\) and  \(i=1,\ldots,n,\)
$$|L_{ni}|\leq \frac{\beta}{\sigma \sqrt{n}},$$ while
$$\sum_{i=1}^n \bE(L_{ni}^2|\calF_{n\,i-1})=1$$ with probability one.

Let  $S_n =\sum_{i=1}^n L_{ni}.$ Then, by the large deviations result \cite[Theorem~1]{Rackauskas90}
for martingales and \eqref{eq:mart_diff_arr} we get
\[
 \limsup_{n\to\infty} \frac{\PP \big(\sum_{i=1}^n {{R_{i}}} \leq nm_\nu -\sigma n^{5/8}\big)}{
\Phi(-n^{1/8})} \leq
 \lim_{n\to\infty} \frac{\PP \big(S_n \leq -n^{1/8}\big)}{
\Phi(-n^{1/8})} =1,
\]
where \(\Phi\) indicates the standard normal distribution.
Since $$\Phi(-x)\leq \frac{\exp(-x^2/2)}{x\sqrt{2\pi}}$$ for all $x>0$,
we obtain
\begin{equation}\label{eq:large_deviation}
 \PP \Big(\sum_{i=1}^n {{R_{i}}}\leq nm_\nu -\sigma n^{5/8}\Big)
\leq \frac{\exp(-n^{1/4}/2)}{n^{1/8}},
\qquad \text{eventually}.
\end{equation}

For \(n=1,2,...,\) set  $F_n=\{\sum_{i=1}^n {{R_{i}}} \leq nm_\nu -\sigma n^{5/8}\} ;$ then
\begin{align*}
 \bE\frac{1}{(c+D_n)^\alpha} & =
 \frac{1}{(c+D_0)^\alpha} \bE\Big(\frac{1}{\big(1+\frac{\sum_{i=1}^n
R_{i}}{c+D_0}\big)^\alpha};F_n\Big)
\\
& \qquad\qquad
+ \bE\Big(\frac{1}{(c+D_0+\sum_{i=1}^n {{R_{i}}})^\alpha};F_n^c\Big)
\\
& \leq
 \frac{\PP (F_n)}{(c+D_0)^\alpha}
+ \frac{1}{(c+D_0+ nm_\nu - \sigma n^{5/8})^\alpha}
\\
& \leq
\frac{1}{(c+D_0 +nm_\nu)^\alpha} \Big( \big(\frac{c+D_0+m_\nu}{c+D_0}\big)^\alpha n^\alpha \PP (F_n)
+ \frac{1}{(1 + o(1))^\alpha} \Big).
\end{align*}
Hence, by \eqref{eq:large_deviation},
\[
 \limsup_{n\to\infty}
\frac{\displaystyle\bE\frac{1}{(c+D_n)^\alpha}}{\displaystyle\frac{1}{(c+D_0+m_Yn)^\alpha}} \leq 1,
\]
and this completes the proof of this case.
\medskip

\noindent
\emph{Case 2: $\sigma=0.$} Assume $\sigma^2_\mu>0$
(the case $\sigma^2_\nu>0$ is analogous).
Hence \(\nu\) is the point mass at \(m_\nu>0.\)  Let $\{\widetilde{R}_Y(n)\}$ be a sequence of independent and identically distributed random variables,
independent of $\{\calA_n\}$ and such that each variable \(\widetilde{R}_Y(n),\) when multiplied by \(m_\mu/m_\nu,\) has probability distribution equal to $\mu.$
For \(n=1,2,...,\) define
\[
\tau^X_n=
 \inf\Big\{k\colon \sum_{i=1}^k \delta_i \geq n\Big\},
\qquad
\tau^Y_n=
 \inf\Big\{k\colon \sum_{i=1}^k (1-\delta_i) \geq n\Big\}.
\]
By Jensen's inequality and \cite[Proposition~2.4]{MayFl07}, we then have:
\begin{align*}
 \bE\frac{1}{(c+D_{n})^\alpha} & =
 \bE\bE\Big(\frac{1}{(c+D_{n})^\alpha}\Big|\sum_{i=1}^n\delta_i=k;\calA_n\Big)
 \\
& =
 \bE\bE\Big(\frac{1}{\big(c+D_0+\sum_{i=1}^k R_X(\tau_i) +
 (n-k)m_\nu\big)^\alpha}\Big|\sum_{i=1}^n\delta_i=k;\calA_n\Big)
\\
& \leq
 \bE\bE\Big(\frac{1}{\big(c+D_0+\sum_{i=1}^k R_X(\tau^X_i) +
 \sum_{i=1}^{n-k} \widetilde{R}_Y(\tau^Y_i)
\big)^\alpha}\Big|\sum_{i=1}^n\delta_i=k;\calA_n\Big)
\\
& =
 \bE\frac{1}{\big(c+D_0+\sum_{i=1}^n (\delta_iR_X(i)+
(1-\delta_i)\widetilde{R}_Y(i)) \big)^\alpha}.
\end{align*}
Since $\min(\sigma^2_\mu,\mathrm{Var}(\widetilde{R}_Y))=(\frac{m_\nu}{m_\mu})^2\sigma^2_\mu>0,$
Case 1 applied to a coupled RRU with the same initial composition and reinforcements equal to \(R_X(n)\) whenever \(\delta(n)=1,\) and  $\widetilde{R}_Y(n) $ whenever \(\delta_n=0,\)
yields the thesis.
\end{pf}

\begin{lemma}\label{subeq:hyp1}
Assume $m_\mu=m_\nu=m>0.$ Then
\[\bE\big(\sum_{k> 0}\sqrt{k} |\Delta A_k|\big)<\infty.\]
\end{lemma}
\begin{pf}
Lemma~\ref{lem:MayGarciaBis} and Lemma~\ref{lem:A^*} yield
\[
\sum_{k> 0} \sqrt{k} \bE(|\Delta A_k|) \leq K_1 \sum_{k> 0} \sqrt{k} \bE\frac{1}{(m+D_k)^2}
\leq K_2 \sum_{k> 0} k^{-3/2} < \infty,
\]
for suitable constants \(K_1,K_2>0.\)
\end{pf}

\begin{lemma}\label{subeq:hyp3}
\[
\bE\big(\sum_{k=0}^\infty k^2 Q_k^4\big) < \infty,\qquad
\bE\big(\sum_{k=0}^\infty k^2 (Q^X_k)^4\big) < \infty,\qquad
\bE\big(\sum_{k=0}^\infty k^2 (Q^Y_k)^4\big) < \infty.
\]
\end{lemma}
\begin{pf}
For all $x>0$ and $0\leq a\leq b$,
\[
\Big(
\frac{a}{b}
\Big)^4 \leq
\Big(
\frac{a+x}{b+x}
\Big)^4,
\]
if $0/0$ is set equal to 1. Then, for \(k=0,1,...,\)
\begin{equation}\label{eq:tobreplaced}
Q_k^4
=
\Big(
\frac{R_{k+1}}{\sum_{i=1}^{k+1}R_i}
\Big)^4
\leq
\Big(
\frac{R_{k+1}+1+D_0}{1+D_{k+1}}
\Big)^4
\leq (1+D_0+\beta)^4
\Big(
\frac{1}{1+D_{k}}
\Big)^4 .
\end{equation}
It follows from Lemma~\ref{lem:MayGarciaBis} with $\alpha=4$ that
\[
\sum_{k=0}^\infty k^2 \bE (Q_k^4) \leq K \sum_{k>0}
k^{-2} <\infty
\]
for a suitable constant \(K>0.\)
Hence $\bE\big(\sum_{k>0}k^2 Q_k^4\big)<\infty.$ The proof is similar for $Q^X$ (resp.\ $Q^Y$): replace
$R_{k+1}$ with $R_X(k+1)$ (resp.\ $R_Y(k+1)$) in the numerator of the first two terms of \eqref{eq:tobreplaced}.
\end{pf}

The next Lemma is an auxiliary result which will be used for proving almost sure convergence of random series.

\begin{lemma}\label{lem:asymp_1.1}
Let \(\{a_k\},\{b_k\}\) and \(\{c_k\}\) be three infinite sequences of real, nonnegative numbers
such that\ $b_k$ and $c_k$ are eventually strictly positive,
$b_k\sim c_k$ and
$\sum_k a_k/b_k < \infty$. Then,
\[
\sum_{k>n} \frac{a_k}{b_k} \sim \sum_{k>n} \frac{a_k}{c_k}
\qquad \text{as }n\to\infty\,.
\]
\end{lemma}
\begin{pf} For lack of a reference, we prove the lemma.
For a fixed $0<\epsilon\leq 1/2$, let $n_0$ be large enough that
$b_k$ and $c_k$ are strictly positive and
$(1-\epsilon)b_k\leq c_k \leq(1+\epsilon)b_k,$ for \(k >n_0.\)
Then, for $n\geq n_0,$
\[
(1-2\epsilon)\sum_{k>n} \frac{a_k}{b_k} \leq
\sum_{k>n} \frac{a_k}{(1+\epsilon)b_k} \leq \sum_{k>n} \frac{a_k}{c_k}
\leq \sum_{k>n} \frac{a_k}{(1-\epsilon)b_k} \leq
(1+2\epsilon)\sum_{k>n} \frac{a_k}{b_k} \,.
\]
\end{pf}

Finally, we need a general fact about convergence of random sequences; for lack of a better reference, see \cite[Lemma 3.2]{PemantleVolkov}.

\begin{lemma}\label{lem:conv_sum}
Let $\{\xi_n\}$ be a sequence of real random variables adapted to the filtration $\{\calA_n\}$.
If \(\PP(\xi_1<\infty)=1\) and
$$\sum_n \bE({{\xi}}_{n+1}|{\mathcal{A}}_n) <\infty \quad \mbox{and}\quad \sum_n \bE({{\xi}}_{n+1}^2|{\mathcal{A}}_n) <\infty $$
almost surely,
then $\sum_n {\xi}_n$ converges almost surely.
\end{lemma}

We can now demonstrate a proposition that will act as cornerstone for the proof of the main result of the paper.

\begin{proposition}\label{subeq:hyp2}
Assume $m_\mu=m_\nu=m>0$, and let
\[
H_X = m^{-2}\bE(R_X^2),\qquad H_Y = m^{-2}{\bE(R_Y^2)}.
\]
Then
\[
\lim_{n\to\infty} n\sum_{k> n} (Q^X_k)^2 = H_X, \qquad \lim_{n\to\infty} n\sum_{k> n} (Q^Y_k)^2 = H_Y
  \]
on a set of probability one.
\end{proposition}

\begin{pf}
We prove that $\lim_{n\to\infty} n\sum_{k> n} (Q^X_k)^2 = H_X$ almost surely, along the argument used to prove
Corollary 4.1 in \cite{Crimaldi08}. The proof that $\lim_{n\to\infty} n\sum_{k> n} (Q^Y_k)^2 = H_Y$ almost surely is similar and
will be omitted.
Let $\rho = \bE(R^2_X)$.
The series
\[
\sum_n n^{-1} (R^2_X(n+1)-\rho)
\]
converges almost surely, since it is a series of zero-mean independent random variables with variances bounded by $n^{-2}\beta^4$.
This fact and Abel's Theorem imply that
\[
 \lim_{n \rightarrow \infty} n\sum_{k {{>}} n} k^{-2} (R^2_X(k+1)-\rho) = 0, 
\]
on a set of probability one.
Then
\begin{equation}\label{eq:Abel}
\lim_{n \rightarrow \infty} n\sum_{k {{>}} n} k^{-2} R^2_X(k+1)= \rho 
\end{equation}
on a set of probability one, since $\lim_{n \rightarrow \infty} n\sum_{k > n} k^{-2} =1.$

From \cite[Lemma~A.1(iii)]{MayFl07}, it follows that $\lim_{k \rightarrow \infty} (mk)^{-1}\sum_{i=1}^k R_i=1$ almost surely and thus
$$(Q^X_k)^2 \thicksim  m^{-2}k^{-2}R^2_X(k+1), 
$$
on a set of probability one.
Therefore Lemma \ref{lem:asymp_1.1} implies that
$$n \sum_{k \ge n} (Q^X_k)^2 \thicksim m^{-2}n \sum_{k \ge n}k^{-2}R^2_X(k+1) 
$$
almost surely; however \eqref{eq:Abel} shows that the right term converges almost surely to $m^{-2} \rho=H_X$ as $n\rightarrow \infty.$ This concludes the proof of the proposition.
\end{pf}
\medskip

\begin{pf}[Proof of Theorem~\ref{thm:mainExtension}]
For \(n=0,1,2,...,\) set $$G_n = \sum_{k> n} \sqrt{k}|\Delta A_k|\ge 0$$ and $W_n=\bE(G_n|\calA_n).$
Because of Lemma~\ref{subeq:hyp1}, the process $\{G_n\}$ converges monotonically to  zero almost surely and in
$L^1,$ as $n$ goes to infinity. Hence the process $\{W_n\}$, being a non-negative super-martingale,
converges to  zero almost surely and in $L^1,$ as $n$ goes to infinity.
Since, for \(n=1,2,...,\)
$$\bE(\sqrt{n}|A_\infty-A_n||{\mathcal A}_n)\leq
\bE(\sum_{k> n} \sqrt{k}|\Delta A_k||{\mathcal A}_n)=\bE(G_n|\calA_n),$$
almost surely, we obtain that, for all $t>0$,
\[
\PP(\sqrt{n}|A_\infty-A_n|>t|{\mathcal A}_n) \leq \frac{\bE(\sqrt{n}|A_\infty-A_n||{\mathcal A}_n)}{t}
\leq \frac{\bE(G_n|\calA_n)}{t}
\]
on a set of probability one;
therefore ${\mathcal{L}}(\sqrt{n}|A_\infty-A_n||{\mathcal A}_n)(\omega)$ weakly converges to the mass function at $0$, for almost every \(\omega \in \Omega.\)
Proving the theorem is thus equivalent to show that, for almost every \(\omega \in \Omega,\)
\({\mathcal{L}} (\sqrt{n}(M_n-M_\infty)|{\mathcal A}_n)(\omega)\)
weakly converges to a \({\mathcal N}(0,\overline{H}(\omega)),\)
where
\[\overline{H}(\omega)=H(\omega)Z_\infty(\omega)(1-Z_\infty(\omega)).\]
Since $\{M_n\}$ is a martingale, this follows from
\cite[Proposition~2.2]{Crimaldi08} once we show that
\begin{subequations}
\begin{align}
&\bE\big(\sup_{k}\sqrt{k} |\Delta M_k|) < \infty \label{subeq:CriP1}
\\
\intertext{and}
&\lim_{n \rightarrow \infty} n\sum_{k> n} (\Delta M_k)^2
  =\overline{H}\;\; \mbox{almost surely}. \label{subeq:CriP2}
\end{align}
\end{subequations}

\noindent\emph{Proof of \eqref{subeq:CriP1}}.
Since
\[
 \sqrt{k}|\Delta M_k| \leq \sqrt{k}|\Delta A_k| +\sqrt{k}|\Delta Z_k|
\]
and
\[
 \bE(\sup_k \sqrt{k}|\Delta A_k|)
\leq \sum_k \sqrt{k}\bE(|\Delta A_k|),
\]
from Lemma~\ref{subeq:hyp1} we get that
\[
\bE\big(\sup_{k}\sqrt{k} |\Delta M_k|)<\infty \iff \bE\big(\sup_{k}\sqrt{k}|\Delta Z_k|)<\infty.
\]
Note that, for \(n=0,1,2,...,\)  $\delta_{n+1}R_X(n+1)=\delta_{n+1}{{R_{n+1}}}$ and
\begin{equation}\label{eq:9est}
\begin{aligned}
Z_n-Z_{n+1} &= \frac{X_n}{D_{n}}-\frac{X_{n+1}}{D_{n+1}}\\
& = \frac{1}{D_nD_{n+1}} \Big( X_nD_{n+1} - X_{n+1}D_{n}\Big)
\\
& = \frac{1}{D_nD_{n+1}} \Big( X_n(D_n+{{R_{n+1}}}) -
(X_{n}+\delta_{n+1}R_X(n+1))D_{n}\Big)
\\
& = \frac{1}{D_nD_{n+1}} \Big( X_n {{R_{n+1}}} -
\delta_{n+1}R_X(n+1)D_{n}\Big)
\\
& = \frac{1}{D_nD_{n+1}} \Big( X_n {{R_{n+1}}} -
\delta_{n+1}{{R_{n+1}}}D_{n}\Big)
\\
& = \frac{{{R_{n+1}}}}{D_{n+1}} \Big( Z_n -\delta_{n+1}\Big)
\\
& = Q_n \frac{\sum_{i=1}^{n+1} R_i}{D_{n+1}}\Big( Z_n -\delta_{n+1}\Big)
\end{aligned}
\end{equation}
which yields
\(|\Delta Z_n|\leq Q_n.\) Hence
\(\bE\big(\sup_{k}\sqrt{k}|\Delta Z_k| \big)^4\leq \bE\big(\sum_k k^2 Q_k^4)<\infty \)
by Proposition~\ref{subeq:hyp2}.
Since $\big(\bE\sup_{k}\sqrt{k}|\Delta Z_k|\big)^4\leq \bE\big(\sup_{k}\sqrt{k}|\Delta Z_k| \big)^4<\infty$ this proves \eqref{subeq:CriP1}.

\null

\noindent\emph{Proof of \eqref{subeq:CriP2}}. We split the proof in four steps.

\noindent
\emph{First step}: We show that
\[
\lim_{n \rightarrow \infty} n\sum_{k> n} (\Delta M_k)^2= \overline{H}\; \mbox{almost surely}
  \iff
\lim_{n \rightarrow \infty}n\sum_{k> n} (\Delta Z_k)^2 =\overline{H}\; \mbox{almost surely}.
\]

Lemma~\ref{subeq:hyp1} shows that
$\bE\big(\sum_{k> 0}\sqrt{k} |\Delta A_k|\big) <\infty$ almost surely; hence
$$\bE\big(\sum_{k> 0}k |\Delta A_k|^2 \big) <\infty$$ almost surely and this implies that
\(\lim_{n \rightarrow \infty} \sum_{k> n}k |\Delta A_k|^2 =0\) on a set of probability one.
However
\(
n |\Delta A_k|^2 \leq k |\Delta A_k|^2,\) for $k> n=1,2,...,$ and thus
\begin{equation}\label{eq:nA^2_kto0}
\lim_{n \rightarrow \infty} n\sum_{k> n} |\Delta A_k|^2=0 \;\;\mbox{almost surely.}
\end{equation}

For $n=1,2,...,$ $\sum_{k\geq 0}k^2 Q_k^4 \geq (\sqrt{n}(\sup_{k> n} Q_k ))^4,$ and
thus Lemma~\ref{subeq:hyp3} implies that $\PP(\sup_n \sqrt{n}(\sup_{k> n} Q_k )=\infty)=0$
which in turn implies, through equation \eqref{eq:9est}, that
$\PP(\sup_n \sqrt{n}(\sup_{k> n} \Delta Z_k)=\infty)=0.$
Hence
\begin{equation}\label{eq:nA^2_kto01}
\lim_{n \rightarrow \infty} \Big|n \sum_{k> n} \Delta Z_k \Delta A_k \Big|\leq \big|\sup_n\sqrt{n}\sup_{k> n}
\Delta Z_k
\Big| \lim_{n \rightarrow \infty} \sqrt{n}\sum_{k> n} |\Delta A_k|
= 0
\end{equation}
almost surely, where the last equality follows, once again, from Lemma~\ref{subeq:hyp1}.
Since, for \(n=1,2,...,\)
\begin{equation}\label{eq:DMZA}
(\Delta M_n)^2 = (\Delta Z_n-\Delta A_n)^2 = (\Delta Z_n)^2 + (\Delta A_n)^2 -2 \Delta Z_n \Delta A_n,
\end{equation}
\eqref{eq:nA^2_kto0} and \eqref{eq:nA^2_kto01} imply that
\[
\lim_{n \rightarrow \infty} n \sum_{k> n} ((\Delta M_n)^2 -(\Delta Z_n)^2 ) = \lim_{n \rightarrow \infty} n \sum_{k> n} ((\Delta A_n)^2 -2 \Delta Z_n \Delta A_n) = 0
\]
on a set of probability one. This concludes the proof of the first step.

For the next three steps, we follow the arguments in \cite[Theorem 1.1]{Crimaldi08}
armed with the results provided by Proposition~\ref{subeq:hyp2} and Lemma~\ref{subeq:hyp3}.

\null
\noindent\emph{Second step}: We show that
\[
\lim_{n \rightarrow \infty}n\sum_{k> n} (\Delta Z_k)^2
  =\overline{H}\;\mbox{almost surely}
  \iff
\lim_{n \rightarrow \infty}n\sum_{k> n} (Z_k-\delta_{k+1})^2
Q^2_k
=\overline{H}\; \mbox{almost surely}.
\]
Lemma~\ref{lem:conv_sum} and \eqref{eq:DMZA} imply the almost sure convergence of
\(\sum_{n} (\Delta Z_{n})^2.\)
Thus, from \eqref{eq:9est} and Lemma~\ref{lem:asymp_1.1}, we get that
\[
\sum_{k> n} (\Delta Z_{k+1})^2 =
\sum_{k> n} (Z_k-\delta_{k+1})^2
\frac{{{R^2_{k+1}}}}{D^2_{k+1}} \sim \sum_{k> n} (Z_k-\delta_{k+1})^2
Q^2_k \qquad 
\]
as \(n\) grows to infinity; this completes the proof of the second step.

\null

\noindent\emph{Third step}: We show that the almost sure convergence of
\begin{equation}\label{eq:serie1}
\sum_{k=0}^m k(\delta_{k+1}-Z_{k}) (1-Z_k)^2 (Q^X_k)^2
\end{equation}
and
\begin{equation}\label{eq:serie2}
\sum_{k=0}^m k(\delta_{k+1}-Z_{k}) Z_k^2 (Q^Y_k)^2,
\end{equation}
as $m$ grows to infinity, implies
that $\lim_{n\to\infty} n\sum_{k> n} (Z_k-\delta_{k+1})^2
Q^2_k =\overline{H}$ almost surely.

Because of Abel's Theorem, almost sure convergence of the series
(\ref{eq:serie1}) and (\ref{eq:serie2}) implies that
\begin{equation}\label{eq:Cri13}
\begin{gathered}
\lim_{n\to\infty} n\sum_{k> n} (\delta_{k+1}-Z_{k}) (1-Z_k)^2 (Q^X_k)^2
= 0,\\
\lim_{n\to\infty} n\sum_{k> n} (\delta_{k+1}-Z_{k}) Z_k^2 (Q^Y_k)^2
= 0
\end{gathered}
\end{equation}
on a set of probability one. Now, from Proposition~\ref{subeq:hyp2} and the almost sure convergence
of the sequence $\{Z_n\}$ to $Z_\infty$, we obtain that,
\begin{equation}\label{eq:Cri12}
\begin{gathered}
\lim_{n\to\infty}n\sum_{k> n} Z_k (1-Z_k)^2
(Q^X_k)^2 = H_X Z_\infty (1-Z_\infty)^2
,\\
\lim_{n\to\infty}n\sum_{k> n} (1-Z_k) Z_k^2
(Q^Y_k)^2 = H_Y (1-Z_\infty) Z_\infty^2.
\end{gathered}
\end{equation}
on a set of probability one.
Equations~\eqref{eq:Cri13}-\eqref{eq:Cri12} yield
\begin{equation*}
\begin{gathered}
\lim_{n\to\infty}n\sum_{k> n} \delta_{k+1} (1-Z_k)^2
(Q^X_k)^2 = H_X Z_\infty (1-Z_\infty)^2
\\
\lim_{n\to\infty} n\sum_{k> n} (1-\delta_{k+1}) Z_k^2
(Q^Y_k)^2 = H_Y (1-Z_\infty) Z_\infty^2
\end{gathered}
\end{equation*}
almost surely.
Since, for all \(k\geq 0,\)
$Q_{k} =\delta_{k+1} Q_{k}^X+(1-\delta_{k+1})Q_{k}^Y$ and  $\delta_{k+1}(1-\delta_{k+1})=0,$ we have
\begin{align*}
\lim_{n\to\infty} n\sum_{k> n} & (Z_k-\delta_{k+1})^2 Q^2_k \\
& = \lim_{n\to\infty} n\sum_{k> n} \Big(
\delta_{k+1} (1-Z_k^2) (Q^X_k)^2 +
(1-\delta_{k+1}) Z_k^2 (Q^Y_k)^2\Big)
\\
& = H_X Z_\infty (1-Z_\infty)^2 +
H_Y (1-Z_\infty) Z_\infty^2 = \overline{H}
\end{align*}
on a set of probability one.

\null
\noindent\emph{Fourth step}:
We prove the almost sure convergence of the series
\[
\sum_{k=0}^{\infty} k(\delta_{k+1}-Z_{k}) (1-Z_k)^2 (Q^X_k)^2;
\]
the proof of the almost sure convergence of $\sum_{k=0}^\infty k(\delta_{k+1}-Z_{k}) Z_k^2 (Q^Y_k)^2$
is similar.

For \(n=0,1,2,...,\) $R_X(n+1)$ is independent of $\sigma(\delta_{n+1},{\mathcal{A}}_n)$ and thus
\[
\bE( n(\delta_{n+1}-Z_{n}) (1-Z_n)^2 \frac{R^2_X({n+1})}{(D_{n}-D_0)^2}|{\mathcal{A}}_n)=0.
\]
Hence,
\begin{eqnarray*}
\lefteqn{\big|\bE(n(\delta_{n+1}-Z_{n}) (1-Z_n)^2 (Q^X_n)^2|{\mathcal{A}}_n)\big|}\\
&&=\Big|\bE(n(\delta_{n+1}-Z_{n}) (1-Z_n)^2 \Big((Q^X_n)^2-\frac{R^2_X({n+1})}{(D_{n}-D_0)^2}\Big) |{\mathcal{A}}_n)\Big|\\
&&\leq\bE(n|\delta_{n+1}-Z_{n}| (1-Z_n)^2 R^2_X({n+1}) \Big(\frac{1}{(D_{n}-D_0)^2}-\frac{1}{(D_{n+1}-D_0)^2}\Big) |{\mathcal{A}}_n)\\
&&\leq n\beta^2 \bE(\frac{1}{(D_{n}-D_0)^2}-\frac{1}{(D_{n+1}-D_0)^2} |{\mathcal{A}}_n)\\
&&\leq 2n\beta^3 \frac{1}{(D_{n}-D_0)^3};
\end{eqnarray*}
the last inequality holds because,
\begin{align*}
\Big(\frac{1}{D_{n}-D_0}\Big)^2-\Big(\frac{1}{D_{n+1}-D_0}\Big)^2
&= \frac{(D_{n+1}-D_0)^2-(D_{n}-D_0)^2}{(D_{n}-D_0)^2(D_{n+1}-D_0)^2}
\\
&\leq \frac{2(D_{n+1}-D_0)R_{n+1}}{(D_{n}-D_0)^2(D_{n+1}-D_0)^2}
\\
&\leq \frac{2\beta}{(D_{n}-D_0)^3}.
\end{align*}
However $\lim_{n\to\infty} D_n/n =m$ almost surely, as proved in Lemma~\cite[Lemma~A.1(iii)]{MayFl07}; thus
$ \sum_n n\bE((\delta_{n+1}-Z_{n}) (1-Z_n)^2 (Q^X_n)^2|{\mathcal{A}}_n) <\infty$ on a set of probability one.

Next note that, as in \cite[Eq.~(16)]{Crimaldi08},
\begin{align*}
\bE\Big( \sum_n n^2\bE((\delta_{n+1}-Z_{n})^2 (1-Z_n)^4 (Q^X_n)^4|{\mathcal{A}}_n) \Big)
\leq \bE\Big( \sum_n n^2 (Q^X_n)^4 \Big) <\infty
\end{align*}
because of Lemma~\ref{subeq:hyp3}. Therefore Lemma~\ref{lem:conv_sum} implies that the series (\ref{eq:serie1})
converges on a set of
probability one; this concludes the proof of the fourth step and that of the theorem.
\end{pf}

\begin{pf}[Proof of Theorem \ref{theo:nomasses}]
Recall that, if $\pi$ and \(\pi'\) are probability
distribution on $\mathbb R$
the discrepancy metric $d_D$ between $\pi'$ and $\pi$ is
defined as
\[
d_D (\pi',\pi) = \sup_{\text{closed balls B}} |\pi'(B)-\pi(B)|;
\]
this metric metrizes the weak convergence of a sequence of probability distributions \(\{\pi_n\}\) to \(\pi,\) when the limiting probability distribution \(\pi\)  is absolutely continuous with respect to Lebesgue measure on $\mathbb R$ (see, e.g., \cite{Gibbs.et.al.02}).

The definition of \(Z_\infty\) and Theorem \ref{thm:mainExtension} imply the existence of
$\Omega'\in \calA$ such that \(\PP(\Omega')=1,\) and, for all  \(\omega \in \Omega',\)
\begin{subequations}\label{okd:st}
\begin{align}
\label{okd:st:01} & \lim_{n \rightarrow \infty}Z_n(\omega)
=Z_{\infty}(\omega)\\
\intertext{and}
\label{okd:st:02} &\lim_{n \rightarrow \infty}d_D((K_n)(\omega),
{\mathcal N} (0 ,\overline{H}(\omega)))=0.
\end{align}
\end{subequations}

By way of contradiction, assume there is a $p\in [0,1]$ such that  $\PP(Z_{\infty}=p)>0$.
Since
$$
\lim_{n \rightarrow \infty}\PP(Z_{\infty}=p|\calA_n)
=1_{\{p\}}(Z_\infty)
$$
almost surely, there is a set $F\in\calA$,
$F\subseteq \{Z_\infty=p\} \bigcap \Omega',$ such that  $\PP(F)>0$ and,
for all $\omega\in F$,
\begin{equation}
\label{nomass:1}
\lim_{n \rightarrow \infty}\PP(Z_{\infty}=p|\calA_n)(\omega)
=1.
\end{equation}

Fix $\omega\in F.$ For \(n=1,2,...,\) set $x_n= \sqrt{n} (Z_n(\omega)-p)$
and consider the closed ball $B_n=\{x_n\}.$ Then, for \(n=1,2,...,\)
\[
d_D(K_n(\omega),{\mathcal N} (0 ,\overline{H}(\omega)) \geq
|K_n(\omega)(B_n) - {\mathcal N} (0 ,\overline{H}(\omega))(B_n)|
= K_n(\omega)(B_n);
\]
however \(\lim_{n \rightarrow \infty}K_n(\omega)(B_n) = 1,\) because of (\ref{nomass:1}), and this contradicts (\ref{okd:st:02}).
\end{pf}

\begin{remark} The same argument works also to show that the distribution of the limit composition \(V\) of the generalized Polya urns treated in \cite{Crimaldi08} has no point masses, whenever the conditions of Theorem 1.1 in \cite{Crimaldi08} are satisfied.
\end{remark}


\begin{pf}[Proof of Theorem \ref{cor:extension}]
Assume $m_\nu>0,$ otherwise it is trivial to prove that \(\PP(Z_\infty=1)=1.\)

We will work through a coupling argument that considers two randomly reinforced urns with the same initial composition \((x,y).\)
Compositions of the first urn are described by the process \(\{(X_n,Y_n)\}\) defined in (\ref{def1}); the composition process
\(\{(\widetilde{X}_n,\widetilde{Y}_n)\}\) of the second urn is defined by
\begin{equation}\label{def11}
\left \{ \begin{array}{lll}
  \widetilde{X}_{n+1}&=&\widetilde{X}_n+R_X(n+1)\widetilde{\delta}_{n+1},\\
  \widetilde{Y}_{n+1}&=&\widetilde{Y}_n+\widetilde{R}_Y(n+1)(1-\widetilde{\delta}_{n+1}),
\end{array} \right.
\end{equation}
where, for \(n=0,1,2,\ldots,\) $\widetilde{R}_Y(n+1)=R_Y(n+1)+(m_\mu-m_\nu)$ and
\(\widetilde{\delta}_{n+1}\) is the indicator of the event
  \(\{U_{n+1}\leq
  \widetilde{X}_{n}(\widetilde{X}_{n}+\widetilde{Y}_{n})^{-1}\}\).
The two urns are coupled because the random sequences \(\{U_n\}\) and \(\{(V_n,W_n)\}\) defining their dynamics
through equations (\ref{def1}) and (\ref{def11}) are the same.

Note that
$m_\mu=\bE(R_X(n+1))=\bE(\widetilde{R}_Y(n+1));$ hence Theorem~\ref{theo:nomasses} implies that
the distribution of \(\widetilde{Z}_\infty\) has no point masses and, in particular,
$$\PP(\widetilde{Z}_\infty=0)=0.$$

By induction on \(n,\) we show that $\widetilde{X}_n\leq{X}_n$
and that $\widetilde{Y}_n\geq{Y}_n$.
For $n=0$ the claim is obvious because the two urns have the same initial composition. Assume the claim to be true for $n$. Then:
\begin{equation}\label{eq:corol_ext1}
Z_n-\widetilde{Z}_{n}=\frac{X_n}{X_n+Y_n}-\frac{\widetilde{X}_n}{\widetilde{X}_n+\widetilde{Y}_n}
=\frac{X_n\widetilde{Y}_n-\widetilde{X}_nY_n}{(X_n+Y_n)(\widetilde{X}_n+\widetilde{Y}_n)}
\geq 0\,,
\end{equation}
which implies \(\delta_{n+1}\geq\widetilde{\delta}_{n+1}\).
Hence
\begin{align*}
  X_{n+1}-\widetilde{X}_{n+1}&=(X_n-\widetilde{X}_n)+R_X(n+1)(\delta_{n+1}-\widetilde{\delta}_{n+1})
   \geq 0,\\
   \intertext{and}
  \widetilde{Y}_{n+1}-Y_{n+1}&=(\widetilde{Y}_n-Y_{n})+R_Y(n+1)(\delta_{n+1}-\widetilde{\delta}_{n+1})
  \\& \hspace{20mm}\qquad +(m_\mu-m_\nu)(1-\widetilde{\delta}_{n+1}) \geq 0.
\end{align*}
Therefore eq.~\eqref{eq:corol_ext1} holds for all $n$; hence
\(\PP(Z_\infty=0){{\leq}} \PP(\widetilde{Z}_\infty=0)=0\).

What remains to prove is that $\PP(Z_\infty\in(0,1))=0$.
To get this, we can use the same argument as in \cite[Theorem 5]{Muliere.et.al.06}, once
it has been proved that \cite[Eq.~(11) in the Proof of Lemma 4]{Muliere.et.al.06}
holds without the assumption of boundedness away from 0 for the supports of the reinforcement distributions.
Defining $A^*_{n}$ as in equation \eqref{eq:A^*}, this
is tantamount to show that
\begin{equation}
\label{eq:Delta_k}
\lim_{n \rightarrow \infty} \sum_{k=1}^n A^*_{k}
=+\infty\qquad
\end{equation}
on a set of probability one.
However, when $m_\mu>m_\nu$, Lemma~\ref{lem:A^*} shows that, for almost every \(\omega \in \Omega,\) there is \(a>0\) such that
$A^*_n(\omega)\geq a/(D_0 +n\beta)$ eventually; hence \eqref{eq:Delta_k} is true.

\end{pf}

\section{A final remark on absolute continuity}
Having proved that the distribution of the limit proportion \(Z_\infty\) of a RRU has no point masses, when the means of the reinforcement distributions \(\mu\) and \(\nu\) are the same, the next obvious question concerns its absolute continuity with respect to Lebesgue measure.
Theorem \ref{theo:nomasses} implies that $\PP(Z_\infty\in S)=0$,
for all countable sets $S$ in \([0,1].\) The next step would be to show that, if \(S\) is a Lebesgue null set,
then $\PP(Z_\infty\in S)=0.$
Unfortunately, the idea developed in the proof of Theorem \ref{theo:nomasses}
cannot be furtherly exploited to produce such result. In any case,
if the closed balls $B_n$ appearing in the proof are replaced with the ``holes'' of a porous set
(for the link between $\sigma$--porous sets and
measures, see \cite{Mera.et.al.03,Zajivcek05}), it
is possible to show that
\[
\PP(Z_\infty\in S)=0,
\]
for all $\sigma$--porous sets $S$ in \([0,1].\)
Unfortunately this is not enough to prove that the distribution of $Z_\infty$
is absolutely continuous; indeed $T$-measures are
singular with respect to the Lebesgue measure, but they attribute  $0$-measure to any $\sigma$-porous set
(see, e.g., \cite{Prokaj01,Tkadlec86}).



\section{Acknowledgments}
We thank Patrizia Berti, Irene Crimaldi, Luca Pratelli and Pietro
Rigo who spot an oversight appearing in a previous version of this
paper while working on a general Central Limit Theorem for
multicolor generalized Polya urns
\cite{BertiRigoCrimaldiPratelli09}.

\end{document}